\documentclass[10pt,reqno]{amsart}
\usepackage{amssymb,amscd,amsbsy}
\usepackage{amsthm}
\renewcommand{\a}{\alpha}
\renewcommand{\b}{\beta}
\newcommand{\g}{\gamma}

\newcommand{\f}{\varphi}
\renewcommand{\o}{\omega}

\newcommand{\B}{{\mathcal B}}

\newcommand{\F}{{\mathcal F}}

\newcommand{\R}{{\Bbb R}}
\newcommand{\Z}{{\Bbb Z}}

\newcommand{\bS}{{\boldsymbol S}}

\newcommand{\rf}[1]{(\ref{#1})}

\newcommand{\df}{\stackrel{\mathrm{def}}{=}}

\newcommand{\trace}{\operatorname{trace}}

\newcommand{\const}{\operatorname{const}}

\newcommand{\eeq}{\end{equation}}
\newcommand{\beq}{\begin{equation}}
\newcommand{\bay}{\begin{eqnarray}}
\newcommand{\ba}{\begin{align*}}
\newcommand{\ea}{\end{align*}}
\newcommand{\ey}{\end{eqnarray}}
\newcommand{\bey}{\begin{eqnarray*}}
\newcommand{\eey}{\end{eqnarray*}}

\newcommand{\be}{\infty}

\newtheorem{thm}{\hspace{\parindent}Theorem}[section]

\theoremstyle{remark}

\newtheorem*{rem*}{Remark}

\newcommand\dg{\frak D}

\begin{document}

\newcommand{\vse}{\vspace{.2in}}

\title{Almost commuting functions of almost\\ commuting self-adjoint operators}

\maketitle
\begin{center}
\Large
Aleksei Aleksandrov$^{\rm a}$, Vladimir Peller$^{\rm b}$
\end{center}

\begin{center}
\footnotesize
{\it$^{\rm a}$St-Petersburg Branch, Steklov Institute of Mathematics, Fontanka 27, 191023 St-Petersburg, Russia and Department of Mathematics and
Mechanics, Saint Petersburg State University, 28, Universitetski pr., St.Petersburg, 198504, Russia\\
$^{\rm b}$Department of Mathematics, Michigan State University, East Lansing, MI 48824, USA}
\end{center}

\medskip

\newcommand{\mt}{{\mathcal T}}

\footnotesize

{\bf Abstract.} Let $A$ and $B$ be almost commuting (i.e,
$AB-BA\in\bS_1$) self-adjoint operators. We construct a functional calculus
$\f\mapsto\f(A,B)$ for $\f$ in the Besov class $B_{\be,1}^1(\R^2)$. This functional calculus is linear, the operators $\f(A,B)$ and $\psi(A,B)$ almost commute for 
$\f,\,\psi\in B_{\be,1}^1(\R^2)$, $\f(A,B)=u(A)v(B)$ whenever $\f(s,t)=u(s)v(t)$,   and the Helton--Howe trace formula holds. The main tool is triple operator integrals.

\medskip

\begin{center}
{\bf\large Fonctions presque commutant d'op\'erateurs auto-adjoints presque commutant}
\end{center}

\medskip

{\bf R\'esum\'e.} On dit que des op\'erateurs $A$ et $B$ presque commutent si leur commutateur $[A,B]$ appartient \`a classe trace. Pour des op\'erateurs $A$ et $B$  auto-adjoints qui presque commutent nous construisons 
un calcul fonctionnel $\f\mapsto\f(A,B)$, $\f\in B_{\be,1}^1(\R^2)$, o\`u 
$B_{\be,1}^1(\R^2)$ est la classe de Besov. Ce calcul a des propri\'et\'es suivantes: il est lin\'eaire, les op\'erateurs $\f(A,B)$ et $\psi(A,B)$ presque commutent pour toutes les fonctions $\f$ et $\psi$ dans $B_{\be,1}^1(\R^2)$,
$\f(A,B)=u(A)v(B)$ si $\f(s,t)=u(s)v(t)$,
et la formule des traces de Helton et Howe est vraie. L'outil principal est les 
int\'egrales triples op\'eratorielles.

\normalsize

\

\begin{center}
{\bf\large Version fran\c caise abr\'eg\'ee}
\end{center}

\medskip

Soient $A$ et $B$ des op\'erateurs auto-adjoints qui {\it presque commutent}, 
c'est-\`a-dire leur commutateur $[A,B]\df AB-BA$ appartient \`a classe trace $\bS_1$.
Dans \cite{HH} on a obtenu la formule suivante:
$$
\trace\big({\rm i}\big(\f(A,B)\psi(A,B)-\psi(A,B)\f(A,B)\big)\big)=
\iint_{\R^2}\left(\frac{\partial\f}{\partial x}\frac{\partial\psi}{\partial y}-
\frac{\partial\f}{\partial y}\frac{\partial\psi}{\partial x}\right)dP
$$
pour tous les polyn\^omes $\f$ et $\psi$ de deux variables, o\`u $P$ est une mesure sign\'ee boreliienne \`a support compact. Il s'est trouv\'e que la mesure $P$ est absolument continue et
$$
dP(x,y)=\frac1{2\pi}g(x,y)\,dx\,dy,
$$
o\`u $g$ est la fonction principale de Pincus (voir la partie anglaise pour plus d'informations).

Le calcul polynomial $\f\mapsto\f(A,B)$ \'etait \'etendu dans \cite{CP} et
\cite{Pe4}. Dans \cite{Pe4} on a construit le calcul fonctionnel $\f\mapsto\f(A,B)$
pour $\f$ appartenant \`a l'intersection des produit tensoriels projectifs 
${\mathcal C}\df \big(L^\be(\R)\hat\otimes B_{\be,1}^1(\R)\big)\bigcap\big(B_{\be,1}^1(\R)\hat\otimes L^\be(\R)\big)$. Ce calcul est linear, les op\'erateurs $\f(A,B)$ et $\psi(A,B)$ presque commute pour $\f,\,\psi\in{\mathcal C}$. En outre si $\f(s,t)=u(s)v(t)$, alors 
$\f(A,B)=u(A)v(B)$. Finalement, la formule de Helton et Howe si-dessus est vrai pour $\f$ et $\psi$ dans ${\mathcal C}$.

Il \'etait aussi d\'emontr\'e dans \cite{Pe4} que c'est impossible de construire un calcul fonctionnel $\f\mapsto\f(A,B)$ pour toutes les fonctions $\f$ continument 
d\'erivables qui ait les propri\'et\'es si-dessus.

Le but de cette note est d'am\'eliorer les r\'esultats de \cite{Pe4}. 

Pour une fonction $\f$ dans la classe de Besov $B_{\infty,1}^1(\R^2)$ nous 
d\'efinissons l'op\'erateur $\f(A,B)$ comme \c ca:
$$
\f(A,B)=\iint f(x,y)\,dE_A(x)\,dE_B(y),
$$
o\`u $E_A$ et $E_B$ sont les mesures sp\'ectrales des op\'erateurs $A$ et $B$.
La th\'eorie d'int\'egrales doubles op\'eratorielles \'etait d\'evelopp\'ee
par Birman et Solomyak \cite{BS} (voir aussi \cite{Pe2} et \cite{AP}).

Le r\'esultat principal de cette note est bas\'e sur la formule suivante:
\begin{align}
\label{komu}
\big[\f(A,B),Q\big]&=
\iiint\frac{\f(x,y_1)-\f(x,y_2)}{y_1-y_2}\,dE_A(x)\,dE_B(y_1)[B,Q]\,dE_B(y_2)\nonumber
\\[.2cm]
&+
\iiint\frac{\f(x_1,y)-\f(x_2,y)}{x_1-x_2}\,dE_A(x_1)[A,Q]\,dE_A(x_2)\,dE_B(y),
\end{align}
o\`u $\f\in B_{\be,1}^1(\R^2)$ et $A$ et $B$ sont des op\'erateurs auto-adjoints pour lesquels les commutateurs $[A,Q]$ et $[B,Q]$ appartiennent \`a $\bS_1$.

Les int\'egrales triples op\'eratorielles 
\bay
\label{troin}
\iiint\Phi(x_1,x_2,x_3)\,dE_1(x_1)T\,dE_2(x_2)R\,dE_3(x_3)
\ey
\'etaint d\'efinies dans \cite{Pe6} pour les fonctions $\Phi$ dans le produit tensoriel projectif int\'egral; ici $T$ et $R$ sont des op\'erateurs born\'es et $E_1$, $E_2$, $E_3$ sont des mesures sp\'ectrales.
Puis dans \cite{JTT} les int\'egrales triples op\'eratorielles \'etaint d\'efinies pour les fonctions qui appartiennent au produit tensoriel de Haagerup
$L^\be(E_1)\otimes_{\rm h}\!L^\be(E_2)\otimes_{\rm h}\!L^\be(E_2)$ (voir \cite{Pis}). Il \'etait \'etabli dans \cite{ANP} que les conditions $T\in\bS_1$ ou $R\in\bS_1$ et
$\Phi\in L^\be(E_1)\otimes_{\rm h}\!L^\be(E_2)\otimes_{\rm h}\!L^\be(E_2)$ 
n'impliquent pas que l'op\'erateur \rf{troin} appartienne \`a classe trace.

On a d\'efini dans \cite{ANP} les produit tensoriels du type de Haagerup $L^\be(E_1)\otimes_{\rm h}\!L^\be(E_2)\otimes^{\rm h}\!L^\be(E_3)$ et $L^\be(E_1)\otimes^{\rm h}\!L^\be(E_2)\otimes_{\rm h}\!L^\be(E_3)$ (voir la partie anglaise). On a d\'emontr\'e dans \cite{ANP} que si $T\in\bS_1$, $R$ est un op\'erateur born\'e
et $\Phi\in L^\be(E_1)\otimes_{\rm h}\!L^\be(E_2)\otimes^{\rm h}\!L^\be(E_3)$, alors
l'op\'erateur \rf{troin} appartient \`a $\bS_1$. De m\^eme, si $R\in\bS_1$, $T$ est un op\'erateur born\'e
et $\Phi\in L^\be(E_1)\otimes^{\rm h}\!L^\be(E_2)\otimes_{\rm h}\!L^\be(E_3)$, alors
l'op\'erateur \rf{troin} appartient \`a $\bS_1$.

Si $\f\in B_{\be,1}^1(\R^2)$, alors la fonction 
$(x_1,x_2,y)\mapsto\big(\f(x_1,y)-\f(x_2,y)\big)(x_1-x_2)^{-1}$ appartient \`a
l'espace $L^\be(E_A)\otimes_{\rm h}\!L^\be(E_A)\otimes^{\rm h}\!L^\be(E_B)$ et la fonction $(x,y_1,y_2)\mapsto\big(\f(x,y_1)-\f(x,y_2)\big)(y_1-y_2)^{-1}$ appartient \`a
$L^\be(E_A)\otimes_{\rm h}\!L^\be(E_B)\otimes^{\rm h}\!L^\be(E_B)$, voir \cite{ANP}.
Donc les int\'egrales dans \rf{komu} sont bien d\'efinis est appartiennent \`a $\bS_1$.

La formule \rf{komu} implique que si $\f$ et $\psi$ appartiennent \`a $B_{\be,1}^1(\R^2)$ et $A$ et $B$ sont des op\'erateurs auto-adjoints qui presque commutent, alors
\begin{align*}
\big[\f(A,B),\psi(A,B)\big]&=
\iiint
\frac{\f(x,y_1)-\f(x,y_2)}{y_1-y_2}\,dE_A(x)\,dE_B(y_1)[B,\psi(A,B)]\,dE_B(y_2)\nonumber
\\[.2cm]
&+
\iiint
\frac{\f(x_1,y)-\f(x_2,y)}{x_1-x_2}\,dE_A(x_1)[A,\psi(A,B)]\,dE_A(x_2)\,dE_B(y)
\end{align*}
et
$$
\big\|[\f(A,B),\psi(A,B)\big]\big\|_{\bS_1}
\le\const\|\f\|_{B_{\be,1}^1(\R^2)}\|\psi\|_{B_{\be,1}^1(\R^2)}
\big\|[A,B]\big\|_{\bS_1}.
$$

\c Ca implique que la formule des traces de Helton et Howe est vrai pour toutes les fonctions $\f$ et $\psi$ dans $B_{\be,1}^1(\R^2)$.

\begin{center}
------------------------------
\end{center}

\setcounter{section}{0}
\section{\bf Introduction}  

\medskip

The spectral theorem allows one for a pair of commuting (bounded) self-adjoint operators $A$ and $B$ to construct a linear and multiplicative functional calculus 
$$
\f\mapsto\f(A,B)=\int_{\R^2}\f(x,y)\,dE_{A,B}(x,y)
$$
for the class of bounded Borel functions on the plane $\R^2$. Here $E_{A,B}$
is the joint spectral measure of the pair$(A,B)$ defined on the Borel subsets of $\R^2$.

If $A$ and $B$ are noncommuting self-adjoint operators, we can define functions of $A$ and $B$ in terms of double operator integrals
\bay
\label{isch}
\f(A,B)\df\iint_{\R^2}\f(x,y)\,dE_A(x)\,dE_B(y)
\ey
for functions $\f$ that are Schur multipliers with respect to the spectral measure $E_A$ and $E_B$ of the operators $A$ and $B$. The theory of double operator integrals was developed by Birman and Solomyak \cite{BS} (we also refer the reader to \cite{Pe2} and \cite{AP} for double operator integrals and Schur multipliers). It was observed in \cite{ANP}  the Besov space $B_{\infty,1}^1(\R^2)$ 
of functions on $\R^2$ is contained in the space of Schur multipliers with respect to compactly supported spectral measures on $\R$, and so for 
$\f\in B_{\infty,1}^1(\R^2)$, the operator $\f(A,B)$ is well defined by \rf{isch} for bounded self-adjoint operators $A$ and $B$
(see \cite{Pee} for the definition and properties of Besov classes).

In this paper we deal with almost commuting self-adjoint operators. {\it Self-adjoint operators $A$ and $B$ are called almost commuting} if their {\it commutator} 
$[A,B]\df AB-BA$ belongs to trace class. 

In \cite{HH} the following trace formula was obtained for bounded almost commuting self-adjoint operators $A$ and $B$:
\bay
\label{HeHo}
\trace\big({\rm i}\big(\f(A,B)\psi(A,B)-\psi(A,B)\f(A,B)\big)\big)=
\iint_{\R^2}\left(\frac{\partial\f}{\partial x}\frac{\partial\psi}{\partial y}-
\frac{\partial\f}{\partial y}\frac{\partial\psi}{\partial x}\right)dP,
\ey
where $P$ is a signed Borel compactly supported measure that corresponds to the pair $(A,B)$. The formula holds for polynomials $\f$ and $\psi$. 

It was shown in \cite{Pin2} that the signed measure $P$ is absolutely continuous with respect to planar Lebesgue measure and
$$
dP(x,y)=\frac1{2\pi}g(x,y)\,dx\,dy,
$$
where $g$ is the {\it Pincus principal function}, which was introduced in \cite{Pin1}.

In \cite{CP} the polynomial functional calculus for almost commuting self-adjoint operators was extended to a functional calculus for the class of functions $\f=\F\o$ that are Fourier transforms of complex Borel measures $\o$ on $\R^2$ satisfying
$$
\int_{\R^2}(1+|t|)(1+|s|)\,d|\o|(s,t)<\be,
$$
and the Helton--Howe trace formula \rf{HeHo} was extended to the class of such functions.

The problem of constructing a rich functional calculus, which would extend the functional calculus constructed in \cite{CP} and for which trace formula \rf{HeHo} would still hold was considered in \cite{Pe4}. The problem was to find a big class of functions ${\mathcal C}$ on $\R^2$ and construct a functional calculus $\f\mapsto\f(A,B)$, $\f\in{\mathcal C}$, that has the following properties:

{\it
{\em(i)}~ the functional calculus $\f\mapsto\f(A,B)$, $\f\in{\mathcal C}$, is linear;

{\em(ii)}~ if $\f(s,t)=u(s)v(t)$, then $\f(A,B)=u(A)v(B)$;

{\em(iii)}~ if $\f,\,\psi\in{\mathcal C}$, then  
$\f(A,B)\psi(A,B)-\psi(A,B)\f(A,B)\in\bS_1$;

{\em(iv)}~ formula {\em\rf{HeHo}} holds for arbitrary $\f$ and $\psi$ in ${\mathcal C}$.
}

Note that the right-hand side of \rf{HeHo} makes sense for arbitrary Lipschitz functions $\f$ and $\psi$. However, it was established
in \cite{Pe4}  that a functional calculus satisfying (i) - (iii) cannot be defined for all continuously differentiable functions. This was deduced from the trace class criterion for Hankel operators (see \cite{Pe1} and \cite{Pe5}).

On the other hand, in \cite{Pe4} estimates of \cite{Pe2} and \cite{Pe3}  were used to construct a functional calculus satisfying (i) - (iv)  for the class 
${\mathcal C}=\big(L^\be(\R)\hat\otimes B_{\be,1}^1(\R)\big)\bigcap\big(B_{\be,1}^1(\R)\hat\otimes L^\be(\R)\big)$. Here $\hat\otimes$ stands for projective tensor product and $B_{\be,1}^1(\R)$ is a Besov class.

In this paper we considerably enlarge  the class of functions
$\big(L^\be(\R)\hat\otimes B_{\be,1}^1(\R)\big)\bigcap\big(B_{\be,1}^1(\R)\hat\otimes L^\be(\R)\big)$ and construct a functional calculus satisfying (i) - (iv) for the Besov class
$B_{\be,1}^1(\R^2)$ of functions of two variables.


\section{\bf Triple operator integrals}

\medskip

Triple operator integrals
\bay
\label{toi}
\iiint\Phi(x_1,x_2,x_3)\,dE_1(x_1)T\,dE_2(x_2)R\,dE_3(x_3)
\ey
were defined in \cite{Pe6} for spectral measures $E_1$, $E_2$, $E_3$, bounded linear operators $T$ and $R$, and for functions $\Phi$ in the integral projective tensor product $L^\be(E_1)\hat\otimes_{\rm i}L^\be(E_2)\hat\otimes_{\rm i}L^\be(E_3)$.

Later the definition of triple operator integrals was extended
in \cite{JTT} to the class of functions $\Phi$ in the Haagerup tensor product
$L^\be(E_1)\otimes_{\rm h}\!L^\be(E_2)\otimes_{\rm h}\!L^\be(E_3)$ (see \cite{Pis} for imformation about Haagerup tensor products).
However, it was shown in \cite{ANP} that unlike in the case 
$\Phi\in L^\be(E_1)\hat\otimes_{\rm i}L^\be(E_2)\hat\otimes_{\rm i}L^\be(E_3)$, the condition 
that $\Phi$ belongs to $L^\be(E_1)\otimes_{\rm h}\!L^\be(E_2)\otimes_{\rm h}\!L^\be(E_3)$ does not guarantee that if one of the operators $T$ and $R$ is of trace class, then the triple operator integral \rf{toi} belongs to $\bS_1$.

In \cite{ANP} the following Haagerup-like tensor products were introduced:

{\bf Definition.}
{\it A function $\Psi$ is said to belong to the tensor product $L^\be(E_1)\otimes_{\rm h}\!L^\be(E_2)\otimes^{\rm h}\!L^\be(E_3)$ if it admits a representation
\bay
\label{yaH}
\Psi(x_1,x_2,x_3)=\sum_{j,k\ge0}\a_j(x_1)\b_{k}(x_2)\g_{jk}(x_3)
\ey
with $\{\a_j\}_{j\ge0},~\{\b_k\}_{k\ge0}\in L^\be(\ell^2)$ and 
$\{\g_{jk}\}_{j,k\ge0}\in L^\be(\B)$, where $\B$ is the space of bounded operators on $\ell^2$}.

For a bounded linear operator $R$ and 
for a trace class operator $T$,  the triple operator integral
$$
W=\iiint\Psi(x_1,x_2,x_3)\,dE_1(x_1)T\,dE_2(x_2)R\,dE_3(x_3)
$$
was defined in \cite{ANP}
as the following continuous linear functional on the class of compact operators:
\bay
\label{fko}
Q\mapsto
\trace\left(\left(
\iiint\Psi(x_1,x_2,x_3)\,dE_2(x_2)R\,dE_3(x_3)Q\,dE_1(x_1)
\right)T\right)
\ey
and it was shown that
$$
\|W\|_{\bS_1}\le\|\Psi\|_{L^\be\otimes_{\rm h}\!L^\be\otimes^{\rm h}\!L^\be}
\|T\|_{\bS_1}\|R\|,
$$
where $\|\Psi\|_{L^\be\otimes_{\rm h}\!L^\be\otimes^{\rm h}\!L^\be}$ is the infimum of $\|\{\a_j\}_{j\ge0}\|_{L^\be(\ell^2)}\|\{\b_k\}_{k\ge0}\|_{L^\be(\ell^2)}
\|\{\g_{jk}\}_{j,k\ge0}\|_{L^\be(\B)}$
over all representations in \rf{yaH}.

Similarly, the tensor product $L^\be(E_1)\otimes^{\rm h}\!L^\be(E_2)\otimes_{\rm h}\!L^\be(E_3)$ was defined in \cite{ANP} as the class of functions $\Psi$ of the form$$
\Psi(x_1,x_2,x_3)=\sum_{j,k\ge0}\a_{jk}(x_1)\b_{j}(x_2)\g_k(x_3)
$$
where $\{\b_j\}_{j\ge0},~\{\g_k\}_{k\ge0}\in L^\be(\ell^2)$, 
$\{\a_{jk}\}_{j,k\ge0}\in L^\be(\B)$. If $T$ is a bounded linear operator, and $R\in\bS_1$, then the continuous linear functional 
$$
Q\mapsto
\trace\left(\left(
\iiint\Psi(x_1,x_2,x_3)\,dE_3(x_3)Q\,dE_1(x_1)T\,dE_2(x_2)
\right)R\right)
$$
on the class of compact operators determines a trace class operator 
$$
W\df\iiint\Psi(x_1,x_2,x_3)\,dE_1(x_1)T\,dE_2(x_2)R\,dE_3(x_3)
$$
and
$$
\|W\|_{\bS_1}\le
\|\Psi\|_{L^\be\otimes^{\rm h}\!L^\be\otimes_{\rm h}\!L^\be}
\|T\|\cdot\|R\|_{\bS_1}.
$$

\section{\bf Commutators of functions of almost commuting self-adjoint operators}

\medskip

Given a differentiable function $\f$ on $\R^2$, we define the divided differences
$\dg_1f$ and $\dg_2f$ on $\R^3$ by
$$
(\dg_1\f)(x_1,x_2,y)\df\frac{\f(x_1,y)-\f(x_2,y)}{x_1-x_2}\quad\mbox{and}\quad
(\dg_2\f)(x,y_1,y_2)=\frac{\f(x,y_1)-\f(x,y_2)}{y_1-y_2}.
$$
It was shown in \cite{ANP} that if $\f$ is a bounded function on $\R^2$ whose Fourier transform is supported 
in the ball $\{\xi\in\R^2:~\|\xi\|\le1\}$, then
$$
(\dg_1\f)(x_1,x_2,y)=
\sum_{j,k\in\Z}\frac{\sin(x_1-j\pi)}{x_1-j\pi}\cdot\frac{\sin(x_2-k\pi)}{x_2-k\pi}
\cdot\frac{\f(j\pi,y)-\f(k\pi,y)}{j\pi-k\pi}.
$$
$$
\sum_{j\in\Z}\frac{\sin^2(x_1-j\pi)}{(x_1-j\pi)^2}
=\sum_{k\in\Z}\frac{\sin^2(x_2-k\pi)}{(x_2-k\pi)^2}=1,
\quad x_1~x_2\in\R,
$$
and
$$
\sup_{y\in\R}\left\|\left\{\frac{\f(j\pi,y)-\f(k\pi,y)}{j\pi-k\pi}
\right\}_{j,k\in\Z}\right\|_\B\le\const\|f\|_{L^\be(\R)}.
$$

It follows that for $\f\in B_{\be,1}^1(\R^2)$,
$$
\dg_1\f\in L^\be(\R)\otimes_{\rm h}\!L^\be(\R)\otimes^{\rm h}\!L^\be(\R)\quad
\mbox{and}\quad
\|\f\|_{L^\be(\R)\otimes_{\rm h}\!L^\be(\R)\otimes^{\rm h}\!L^\be(\R)}
\le\const\|\f\|_{B_{\be,1}^1}.
$$

Similarly, for $\f\in B_{\be,1}^1(\R^2)$,
$$
\dg_2\f\in L^\be(\R)\otimes^{\rm h}\!L^\be(\R)\otimes_{\rm h}\!L^\be(\R)
\quad\mbox{and}\quad
\|\f\|_{L^\be(\R)\otimes^{\rm h}\!L^\be(\R)\otimes_{\rm h}\!L^\be(\R)}
\le\const\|\f\|_{B_{\be,1}^1}.
$$

\begin{thm}
\label{komu}
Let $A$ and $B$ be self-adjoint operators and let $Q$ be a bounded linear operator such that $[A,Q]\in\bS_1$ and $[B,Q]\in\bS_1$. 
Suppose that $\f\in B_{\be,1}^1(\R^2)$.
Then $[\f(A,B),Q\big]\in\bS_1$,
\begin{align}
\label{komQ}
\big[\f(A,B),Q\big]&=
\iiint\frac{\f(x,y_1)-\f(x,y_2)}{y_1-y_2}\,dE_A(x)\,dE_B(y_1)[B,Q]\,dE_B(y_2)\nonumber
\\[.2cm]
&+
\iiint\frac{\f(x_1,y)-\f(x_2,y)}{x_1-x_2}\,dE_A(x_1)[A,Q]\,dE_A(x_2)\,dE_B(y)
\end{align}
and
$$
\big\|[\f(A,B),Q\big]\big\|_{\bS_1}
\le\const\|\f\|_{B_{\be,1}^1(\R^2)}\big(\big\|[A,Q]\big\|_{\bS_1}+
\big\|[B,Q]\big\|_{\bS_1}\big).
$$
\end{thm}


To obtain the main result of the paper, we apply Theorem \ref{komu} in the case
$Q=\psi(A,B)$, where $\psi\in B^1_{\be,1}(\R^2)$.

\begin{thm}
\label{glav}
Let $A$ and $B$ be almost commuting self-adjoint operators and let $\f$ and $\psi$ be functions in the Besov class $B_{\be,1}^1(\R^2)$. Then
\begin{align}
\label{kom}
\big[\f(A,B),\psi(A,B)\big]&=
\iiint
\frac{\f(x,y_1)-\f(x,y_2)}{y_1-y_2}\,dE_A(x)\,dE_B(y_1)[B,\psi(A,B)]\,dE_B(y_2)\nonumber
\\[.2cm]
&+
\iiint
\frac{\f(x_1,y)-\f(x_2,y)}{x_1-x_2}\,dE_A(x_1)[A,\psi(A,B)]\,dE_A(x_2)\,dE_B(y)
\end{align}
and
\bay
\label{ner}
\big\|[\f(A,B),\psi(A,B)\big]\big\|_{\bS_1}
\le\const\|\f\|_{B_{\be,1}^1(\R^2)}\|\psi\|_{B_{\be,1}^1(\R^2)}
\big\|[A,B]\big\|_{\bS_1}.
\ey
\end{thm}

Note that the right-hand side of inequality \rf{ner} does not involve the norms of $A$ or $B$. Thus formulae \rf{komQ} and \rf{kom} allow us to consider commutators
$\big[f(A,B),g(A,B)\big]$ even for unbounded self-adjoint operators $A$ and $B$
with trace class commutator $[A,B]$, though the the functions $f(A,B)$ and $g(A,B)$ of $A$ and $B$ are not necessarily defined for $f$ and $g$ in 
$B_{\be,1}^1(\R^2)$.

Let us also mention that in \cite{Pe4} it was proved that for almost commuting self-adjoint operators $A$ and $B$, the functional calculus
$\f\mapsto\f(A,B)$, 
$\f\in \big(L^\be(\R)\hat\otimes B_{\be,1}^1(\R)\big)\bigcap\big(B_{\be,1}^1(\R)\hat\otimes L^\be(\R)\big)$, is {\it almost multiplicative}, i.e.,
$$
(\f\psi)(A,B)-\f(A,B)\psi(A,B)\in\bS_1,\quad\f,\psi
\in \big(L^\be(\R)\hat\otimes B_{\be,1}^1(\R)\big)\bigcap\big(B_{\be,1}^1(\R)\hat\otimes L^\be(\R)\big).
$$
It would be interesting to find out whether the functional calculus
$\f\mapsto\f(A,B)$, $\f\in B_{\be,1}^1(\R^2)\bigcap L^\be(\R^2)$, is also almost multiplicative.

\section{\bf An extension of the Helton--Howe trace formula}

\medskip

In this section we use the results of the previous section to extend the Helton--Howe trace formula.

\begin{thm}
Let $A$ and $B$ be almost commuting self-adjoint operators and let $\f$ and $\psi$ be functions in the Besov class $B_{\be,1}^1(\R^2)$. Then the following formula holds:
\bay
\label{exHH}
\trace\big({\rm i}\big(\f(A,B)\psi(A,B)-\psi(A,B)\f(A,B)\big)\big)=\frac{1}{2\pi}
\iint_{\R^2}\left(\frac{\partial\f}{\partial x}\frac{\partial\psi}{\partial y}-
\frac{\partial\f}{\partial y}\frac{\partial\psi}{\partial x}\right)g(x,y)\,dx\,dy,
\ey
where $g$ is the Pincus principal function associated with the operators 
$A$ and $B$.
\end{thm}

It would be interesting to extend the Pincus principal function to the case of unbounded self-adjoint operators with trace class commutators and extend formula \rf{exHH} to unbounded almost commuting operators.

\end{document}